\documentclass[11pt]{amsart}
\usepackage{amssymb,amsmath,txfonts,mathrsfs,mathtools,titletoc, tikz, stmaryrd}
\usepackage[alphabetic]{amsrefs}

\newtheorem{theorem}{Theorem}[section]

\newtheorem{lemma}[theorem]{Lemma}

\newtheorem{question}[theorem]{Question}

\newtheorem{cor}[theorem]{Corollary}

\newtheorem{conj}[theorem]{Conjecture}

\numberwithin{equation}{section}

\def\pf{{\it Proof:}~}

\begin{document}

\title[Integral of scalar curvature]{Integral of scalar curvature on manifolds with a pole}
\author{Guoyi Xu}
\address{Department of Mathematical Sciences\\Tsinghua University, Beijing\\P. R. China, 100084}
\email{guoyixu@tsinghua.edu.cn}
\date{\today}
\date{\today}

\begin{abstract}
On any complete three dimensional Riemannian manifold with a pole and non-negative Ricci curvature,  we show that the asymptotic scaling invariant integral of scalar curvature,  is equal to a term determined by the asymptotic volume ratio of this Riemannian manifold.
\\[3mm]
\end{abstract}
\thanks{The author was partially supported by Beijing Natural Science Foundation Z190003, NSFC 11771230 and NSFC 12141103.}

\maketitle

\section{Introduction}

Since the well-known Gauss-Bonnet Theorem, there are a lot of studies around the integral of scalar curvature.  Especially,  Cohn-Vossen \cite{CV} proved: 
\begin{theorem}\label{thm CV's ineq}
{If $(M^2, g)$ is a finitely connected, complete, oriented Riemannian manifold, and assume $\int_{M^2} R$ exists as an extended real number, then 
\begin{align}
\int_{M^2} R \leq 4\pi\cdot \chi(M^2),  \nonumber 
\end{align}
where $R$ is the scalar curvature and $\chi(M^2)$ is the Euler characteristic of $(M^2, g)$ respectively.
}
\end{theorem}

Motivated to get a generalization of the Cohn-Vossen's inequality, Yau \cite[Problem $9$]{Yau} posed the following question: 
\begin{question}\label{ques Yau}
{Given a $n$-dimensional complete manifold $(M^n, g)$ with $Rc\geq 0$, let $B_p(r)$ be the geodesic ball around $p\in M^n$ and $\sigma_k$ be the $k$-th elementary symmetric function of the Ricci tensor, is it true that $\displaystyle \lim\limits_{r\rightarrow \infty} r^{-n+2k} \int_{B_p(r)} \sigma_k< \infty$?
}
\end{question}
Related to the above question, Shi and Yau \cite{ShiYau} gave a scaling invariant upper bound estimate for the average integral of the scalar curvature (note that the scalar curvature is $\sigma_1$ in Question \ref{ques Yau}), on K\" ahler manifolds with bounded, pinched, nonnegative holomorphic bisectional curvature.  Also see \cite{Petrunin} for discussion with respect to manifolds with lower bound of sectional curvature.

Later Yang \cite{YangBo} constructed examples, which answered Question \ref{ques Yau} for $k> 1$ negatively.  But the case $k= 1$ is left open, where $\sigma_1$ is the scalar curvature $R$. We reformulate the rest case of Question \ref{ques Yau} as follows. 
\begin{question}[Yau]\label{ques Yau-1}
{For any complete Riemannian manifold $(M^n, g)$ with $Rc\geq 0$, any $p\in M^n$, is it true that 
\begin{align}
\lim\limits_{r\rightarrow \infty} r^{2- n}\int_{B_p(r)} R< \infty ? \nonumber 
\end{align}
}
\end{question}

Let $\omega_n$ be the volume of the unit ball in $\mathbb{R}^n$, we recall the definition of the \textbf{asymptotic volume ratio} of the manifold $M^n$ with $Rc\geq 0$: 
\begin{align}
\mathrm{V}_{M^n}\vcentcolon= \lim\limits_{r\rightarrow \infty}\frac{\mathrm{Vol}(B_p(r))}{\omega_nr^n}\in [0, 1]. \label{def of AVR}
\end{align}

Using the monotonicity formulas of Colding and Minicozzi \cite{CM}, we \cite{Xu} obtained: for a complete non-compact, non-parabolic Riemannian manifold $M^3$ with $Rc\geq 0$, there is 
\begin{align}
\varlimsup_{r\rightarrow \infty}\frac{\int_{b\leq r}R\cdot |\nabla b|}{r}\leq 8\pi \big[1- \mathrm{V}_{M^3}\big], \nonumber 
\end{align}
where $b= G^{-1}$ and $G$ is the positive Green's function of $(M^3, g)$.

If there is $p\in (M^n, g)$ such that $\exp_p: T_pM^n\rightarrow M^n$ is a diffeomorphism, we call $(M^n,  g)$ is a Riemannian manifold with a \textbf{pole} (also see \cite{GW}). Among other things, for $3$-dim Riemannian manifolds with a pole and $Rc\geq 0$,  Zhu \cite{Zhu-1}, \cite{Zhu} proved
\begin{align}
\varlimsup_{r\rightarrow \infty} r^{-1}\int_{B_p(r)} R\leq 20\pi. \nonumber 
\end{align}

In this note, we prove the following theorem, which answers Question \ref{ques Yau-1} for $3$-dim manifolds with a pole and $Rc\geq 0$. 
\begin{theorem}\label{thm main result}
{If $(M^3, g)$ is a Riemannian manifold with a pole and $Rc\geq 0$, then 
\begin{align}
\lim_{r\rightarrow\infty} \frac{\int_{B_q(r)}R}{r}= 8\pi (1- \mathrm{V}_{M^3}), \quad \quad \quad \forall q\in M^3. \nonumber 
\end{align}
}
\end{theorem}

Recall that $(M^n, g)$ is called \textbf{Ricci-pinched} if there is $\epsilon> 0$ such that
\begin{align}
Rc\geq \epsilon\cdot R\cdot g\geq 0. \nonumber 
\end{align}

Hamilton has proposed the following conjecture:
\begin{conj}\label{conj Hamilton}
{If $(M^3, g)$ is a complete, non-compact Riemannian manifold and Ricci pinched, then $(M^3, g)$ is flat.
}
\end{conj}
The conjecture has been verified with additional assumption by \cite{CZ} and \cite{Lott}. For recent progress on Conjecture \ref{conj Hamilton}, see \cite{DSS} and \cite{LT} for discussion using Ricci flow (also see \cite{HK} for argument by
the inverse mean curvature flow).

In this note, we prove the following result as a byproduct of the proof of Theorem \ref{thm main result}, which is a special case of Conjecture \ref{conj Hamilton}.
\begin{theorem}\label{thm Hamilton's conj}
{If $(M^3, g)$ is a Riemannian manifold with a pole and Ricci pinched, then $(M^3, g)$ is flat.
}
\end{theorem}

We conclude this section with the sketch of the proof of Theorem \ref{thm main result} and the structure of the paper.  By the Gauss-Codazzi equation,  we firstly decompose the scalar curvature into three parts.  

The first part is the intrinsic scalar curvature of geodesic sphere, whose integral can be calculated by Gauss-Bonnet Theorem;  because geodesic sphere is diffeomorphic to $\mathbb{S}^2$ in our case. 

The second part and third part involve the Ricci curvature along the radial direction and the principal curvatures of the geodesic spheres.  By Calabi's classical computation,  we reduce these two parts as the Ricci curvature along the radial direction and the radial derivative of geodesic spheres' area.  

One difficulty we need to overcome is estimating the average integral of the Ricci curvature along the radial direction.  We firstly observe the derivative of volume elements of geodesic sphere along the radial direction is non-negative on manifolds with a pole.  This observation and Gromov-Bishop's Volume Comparison Theorem, reduce the estimate of radial Ricci curvature to the estimate of the radial derivative of geodesic spheres' area.  

Then we expand (or shrink) the geodesic ball suitably,  such that the radial derivative of suitable geodesic sphere' area can be represented as the difference of two geodesic spheres' area,  by the help of the Mean Value Theorem.  This trick is partially motivated by the argument in \cite{Petrunin}.  Again combining the Gromov-Bishop's Volume Comparison Theorem,  we get the estimate of the radial derivative of geodesic spheres' area.  Furthermore, the average  integral of the Ricci curvature along the radial direction vanishes.

The structure of the paper is as follows.  The estimate of the radial derivative of suitable geodesic sphere' area is discussed in Section \ref{sec volume of sphere}.  Then we present the estimate of Ricci curvature along the radial direction in Section \ref{sec Ricci along the radial direction}.  Finally,  the proof of Theorem \ref{thm main result} and Theorem \ref{thm Hamilton's conj} is provided in Section \ref{sec main proof}.

\section{The volume of geodesic spheres}\label{sec volume of sphere}

Without loss of generality, we assume $p$ is one pole of $(M^3, g)$ in the rest argument.  And we always use $\rho(x)= d(p, x)$ unless otherwise mentioned.  

For any $q\in M^n$,  we can write $q= (r, \theta)$ in terms of polar normal coordinates at $p$, where $r= d(p, q)$.  And it is well known that we can write the volume element of $(M^n,  g)$ as 
\begin{align}
d\mu_g= J(t, \theta)dtd\theta,  \label{measure in polar form}
\end{align}
where $d\theta$ is the area element of the unit $(n- 1)$-sphere. The area element of geodesic sphere $\partial B_p(t)$ is given by $J(t, \theta)d\theta$.

Now we recall the Bishop-Gromov Volume Comparison Theorem for manifolds with a pole (see \cite{LiBook} etc.) as follows:
\begin{theorem}[Bishop-Gromov]\label{thm Bishop-Gromov thm}
{If $(M^n, g)$ is a Riemannian manifold with a pole and $Rc\geq 0$, assume $s\geq t> 0$, then 
\begin{align}
&\frac{J(s, \theta)}{J(t, \theta)}\leq \frac{s^{n- 1}}{t^{n- 1}}, \quad \quad \frac{J'}{J}(t, \theta)\leq \frac{n- 1}{t}, \quad \quad J(t, \theta)\leq t^{n- 1},\nonumber \\
&\frac{V(B_p(s))}{V(B_p(t))}\leq \frac{s^n}{t^n}, \quad \quad \quad V(B_p(t))\leq \omega_n t^n. \nonumber 
\end{align}
}
\end{theorem}\qed

Define $\mathscr{A}(t)= V(\partial B_p(t))$, where $t> 0$. From (\ref{def of AVR}) and Theorem \ref{thm Bishop-Gromov thm},  we get
\begin{align}
\lim_{r\rightarrow\infty}\frac{\mathscr{A}(r)}{4\pi r^2}= \lim_{r\rightarrow\infty}\frac{\int_{\mathbb{S}^2}J(r, \theta)d\theta}{4\pi r^2}= \mathrm{V}_{M^3}. \label{volume limit of level sets}
\end{align}

Fixing $\epsilon\in (0, 1)$, for any $s> 0$, there are $b_s\in ((1- \epsilon)s, s), c_s\in (s, (1+ \epsilon)s)$ such that
\begin{align}
\mathscr{A}'(b_s)= \frac{\mathscr{A}(s)- \mathscr{A}((1- \epsilon)s)}{\epsilon s}, \quad \quad \mathscr{A}'(c_s)= \frac{\mathscr{A}((1+ \epsilon)s)- \mathscr{A}(s)}{\epsilon s}. \nonumber 
\end{align}

For any $s> a$,  we define $A_{a, s}= B_p(s)- B_p(a)$ in the rest argument.

\begin{lemma}\label{lem deri of volume of slice limits}
{If $(M^3, g)$ is a Riemannian manifold with a pole $p$ and $Rc\geq 0$, then for any $a> 0$,  we have 
\begin{align}
\varlimsup_{s\rightarrow\infty}\int_{A_{a, c_s}}\frac{-J''}{s\cdot J}= 4\pi \mathrm{V}_{M^3}(- 2- \epsilon), \quad \quad  \varliminf_{s\rightarrow\infty}\int_{A_{a, b_s}}\frac{-J''}{s\cdot J} =  4\pi \mathrm{V}_{M^3}\cdot (- 2+ \epsilon). \nonumber 
\end{align}
}
\end{lemma}

\pf
{Now using (\ref{volume limit of level sets}), we get 
\begin{align}
&\varlimsup_{s\rightarrow\infty}\frac{1}{s}\int_{A_{a, c_s}}(-\frac{J''}{J})= \varlimsup_{s\rightarrow\infty}\frac{\mathscr{A}'(a)- \mathscr{A}'(c_s)}{s}= \varlimsup_{s\rightarrow\infty}\frac{\mathscr{A}(s)- \mathscr{A}((1+ \epsilon)s)}{\epsilon s^2} \nonumber \\
&= 4\pi \mathrm{V}_{M^3}(- 2- \epsilon). \nonumber \\
&\varliminf_{s\rightarrow\infty}\int_{A_{a, b_s}}\frac{-J''}{s\cdot J} = \varliminf_{s\rightarrow\infty}\frac{\mathscr{A}'(a)- \mathscr{A}'(b_s)}{s} = \varliminf_{s\rightarrow\infty}\frac{\mathscr{A}((1- \epsilon)s)- \mathscr{A}(s)}{\epsilon s^2}. \nonumber \\
&=  4\pi \mathrm{V}_{M^3}\cdot (\epsilon- 2). \nonumber 
\end{align}
}
\qed

\section{The Ricci curvature along the radial direction}\label{sec Ricci along the radial direction}

\begin{lemma}\label{lem bound of deri of J over J}
{If $(M^3, g)$ is a Riemannian manifold with a pole and $Rc\geq 0$, 
\begin{align}
0\leq \frac{J'}{J}(t, \theta)\leq \frac{2}{t}, \quad \quad \quad\quad \forall\ t\in \mathbb{R}^+, \theta\in \mathbb{S}^{n -1}. \nonumber
\end{align}
}
\end{lemma}

\pf
{From \cite{Calabi}, we have 
\begin{align}
\frac{J'}{J}= \lambda_1+ \lambda_2,  \quad \quad \quad (\frac{J'}{J})'= -(\lambda_1^2+ \lambda_2^2)- Rc(\nabla \rho, \nabla \rho),  \nonumber 
\end{align}
where $\lambda_i$ are principal curvature of $\rho^{-1}(t)$ where $i= 1, 2$.  Therefore we get
\begin{align}
Rc(\nabla \rho)=  2\lambda_1\lambda_2- \frac{J''}{J}. \label{upper bound of Rc along ray}
\end{align}

By (\ref{upper bound of Rc along ray}), we have 
\begin{align}
Rc(\nabla \rho)\leq \frac{1}{2}(\lambda_1+ \lambda_2)^2- \frac{J''}{J}= \frac{1}{2}(\frac{J'}{J})^2- \frac{J''}{J}. \label{upper bound of Rc along ray ineq}
\end{align}

From (\ref{upper bound of Rc along ray ineq}), we obtain 
\begin{align}
(-\frac{1}{2}\frac{J'}{J})'\geq \frac{1}{2}Rc(\nabla \rho)+ (\frac{1}{2}\frac{J'}{J})^2\geq (\frac{1}{2}\frac{J'}{J})^2. \label{ODE ineq}
\end{align}
Note $J(t)> 0$ and $J'(t)$ is well-defined for any $t\in \mathbb{R}^+$, because $M^3$ is a Riemannian manifold with a pole. Hence from (\ref{ODE ineq}) we know that 
\begin{align}
J'(t)\geq 0, \quad \quad \quad \quad \forall t\in \mathbb{R}^+. \nonumber 
\end{align}
Otherwise $\frac{-J'}{J}$ will be positive infinity for some finite $t_0> 0$ because of (\ref{ODE ineq}). 

Combining Theorem \ref{thm Bishop-Gromov thm}, we get the conclusion.
}
\qed

\begin{lemma}\label{lem integ of Rc alogn the ray}
{If $(M^3, g)$ is a Riemannian manifold with a pole $p$ and $Rc\geq 0$, then
\begin{align}
\lim_{r\rightarrow\infty} \frac{\int_{B_p(r)}Rc(\nabla \rho)}{r}=0. \nonumber 
\end{align}
}
\end{lemma}

\pf
{\textbf{Step (1)}. By (\ref{volume limit of level sets}), for $\epsilon> 0$, we can find $a> 0$ such that $\displaystyle \frac{\mathscr{A}(a)}{a^2}\leq 4\pi \mathrm{V}_{M^3}(1+ \epsilon)$. By Theorem \ref{thm Bishop-Gromov thm}, we know that $\frac{J(t, \theta)}{t^2}\leq \frac{J(a, \theta)}{a^2}$ for any $t\geq a$. Therefore, integration by parts yields
\begin{align}
&\varlimsup_{s\rightarrow\infty} \frac{1}{s}\int_{A_{a, c_s}} \frac{J'(t, \theta)}{t\cdot J(t, \theta)}= \varlimsup_{s\rightarrow\infty} \frac{1}{s}\int_{\mathbb{S}^2}d\theta\int_a^{c_s} \frac{J'(t, \theta)}{t}dt \nonumber \\
&= \varlimsup_{s\rightarrow\infty} \frac{1}{s}\cdot \Big\{\int_{\mathbb{S}^2}\frac{J(c_s, \theta)}{c_s}d\theta+ \int_{\mathbb{S}^2}d\theta \int_a^{c_s}\frac{J}{t^2}dt\Big\} \nonumber \\
&= \varlimsup_{s\rightarrow\infty} \frac{\mathscr{A}(c_s)}{s\cdot c_s}+ \frac{c_s- a}{s} \frac{\mathscr{A}(a)}{a^2} \nonumber \\
&\leq (1+ \epsilon)\cdot \varlimsup_{s\rightarrow\infty} \frac{\mathscr{A}(c_s)}{c_s^2}+ (1+ \epsilon)\cdot 4\pi \mathrm{V}_{M^3}(1+ \epsilon) \nonumber \\
&\leq 8\pi \mathrm{V}_{M^3}(1+ \epsilon)^2. \label{Rc integ upper bound-1}
\end{align}

\textbf{Step (2)}. Now from (\ref{upper bound of Rc along ray ineq}) and Lemma \ref{lem bound of deri of J over J}, we get
\begin{align}
Rc(\nabla \rho)(t, \theta)\leq \frac{J'(t, \theta)}{t\cdot J(t, \theta)}- \frac{J''}{J}(t, \theta). \label{upper of Rc refined ineq}
\end{align}

Using (\ref{upper of Rc refined ineq}), (\ref{Rc integ upper bound-1}) and Lemma \ref{lem deri of volume of slice limits}, we get 
\begin{align}
&\varlimsup_{s\rightarrow\infty}\frac{1}{s}\int_{B_p(s)}Rc(\nabla \rho)= \varlimsup_{s\rightarrow\infty}\frac{\int_{A_{a, s}}Rc(\nabla \rho)}{s}\leq  \varlimsup_{s\rightarrow\infty} \frac{\int_{A_{a, c_s}} Rc(\nabla \rho)}{s} \nonumber \\
&\leq \varlimsup_{s\rightarrow\infty} \frac{\int_{A_{a, c_s}} \frac{J'(t, \theta)}{t\cdot J(t, \theta)}}{s}+ \varlimsup_{s\rightarrow\infty}\int_{A_{a, c_s}}\frac{-J''}{s\cdot J} \leq 4\pi \mathrm{V}_{M^3}(3\epsilon+ 2\epsilon^2) . \nonumber 
\end{align}

Let $\epsilon\rightarrow 0$ above, we have
\begin{align}
\varlimsup_{s\rightarrow\infty}\frac{1}{s}\int_{B_p(s)}Rc(\nabla \rho)\leq 0. \nonumber 
\end{align}
On the other hand we know that $Rc(\nabla \rho)\geq 0$, therefore the conclusion follows.
}
\qed

\section{The integral of the scalar curvature}\label{sec main proof}

\begin{theorem}
{If $(M^3, g)$ is a Riemannian manifold with a pole and $Rc\geq 0$, then 
\begin{align}
\lim_{r\rightarrow\infty} \frac{\int_{B_q(r)}R}{r}= 8\pi (1- \mathrm{V}_{M^3}), \quad \quad \quad \forall q\in M^3. \nonumber 
\end{align}
}
\end{theorem}

\pf
{\textbf{Step (1)}.  Recall that $\lambda_i$ are principal curvature of $\rho^{-1}(t)$ where $i= 1, 2$.  From the Gauss-Codazzi equation we have 
\begin{align}
R(M^3)(t,  \theta)= R(\partial B_p(t))(t,  \theta)+ 2Rc(M^3)\Big(\nabla \rho, \nabla \rho\Big)- 2\lambda_1\lambda_2,\label{true ineq} 
\end{align}
where $R(\partial B_p(t))(t,  \theta)$ is the intrinsic scalar curvature of $\partial B_p(t)$ at $\exp_p(t\theta)$ and $\exp_p$ is the exponential map of $M^3$.

Note $\rho^{-1}(t)$ is diffeomorphic to $\mathbb{S}^2$.  By the co-area formula and the Gauss-Bonnet Theorem for the compact surface $\rho^{-1}(t)$,  we get 
\begin{align}
\int_{\rho^{-1}[a, s]}R(\partial B_p(t))dx&= \int_a^s dt\int_{\rho^{-1}(t)}R(\partial B_p(t))= 8\pi(s- a).  \label{eq- 1.0}
\end{align}

From (\ref{true ineq}) and (\ref{upper bound of Rc along ray}), we get
\begin{align}
\int_{\rho^{-1}[a, s]}R(M^3)&= \int_{\rho^{-1}[a, s]}\Big\{R(\partial B_p(t))+ 2Rc(M^3)\Big(\nabla \rho, \nabla \rho\Big)- 2\lambda_1\lambda_2\Big\}\nonumber\\
&= 8\pi(s- a)+ \int_{\rho^{-1}[a, s]}\Big\{ Rc(\nabla \rho)- \frac{J''}{J}\Big\}. \label{integ equa pf scalar curv}
\end{align}

\textbf{Step (2)}. From (\ref{integ equa pf scalar curv}) and $Rc\geq 0$, using Lemma \ref{lem deri of volume of slice limits} and $\frac{b_s}{s}\in [1- \epsilon, 1]$, we have 
\begin{align}
&\varliminf_{s\rightarrow\infty}\frac{1}{s}\int_{A_{a, s}}R\geq \varliminf_{s\rightarrow\infty}\frac{1}{s}\int_{A_{a, b_s}}R\geq 8\pi(1- \epsilon)+ \varliminf_{s\rightarrow\infty}\int_{A_{a, b_s}}\frac{-J''}{s\cdot J} \nonumber \\
&= 8\pi(1- \epsilon)+ 4\pi \mathrm{V}_{M^3}\cdot (\epsilon- 2). \nonumber 
\end{align}

Let $\epsilon\rightarrow 0$ in the above, we get 
\begin{align}
\varliminf_{s\rightarrow\infty}\frac{1}{s}\int_{B_p(s)}R= \varliminf_{s\rightarrow\infty}\frac{1}{s}\int_{A_{a, s}}R\geq 8\pi(1- \mathrm{V}_{M^3}). \label{lower limit of global integral}
\end{align}

\textbf{Step (3)}. Using (\ref{integ equa pf scalar curv}), Lemma \ref{lem integ of Rc alogn the ray} and Lemma \ref{lem deri of volume of slice limits}, also note $\frac{c_s}{s}\in [1, 1+ \epsilon]$, we gave 
\begin{align}
&\varlimsup_{s\rightarrow\infty}\frac{1}{s}\int_{B_p(s)}R= \varlimsup_{s\rightarrow\infty}\frac{1}{s}\int_{A_{a, s}}R\leq \varlimsup_{s\rightarrow\infty}\frac{1}{s}\int_{A_{a, c_s}}R \nonumber \\
&\leq 8\pi (1+ \epsilon)+ \varlimsup_{s\rightarrow\infty} \frac{c_s}{s}\cdot \frac{1}{c_s}\int_{B_p(c_s)} Rc(\nabla \rho)+ \varlimsup_{s\rightarrow\infty}\frac{1}{s}\int_{A_{a, c_s}}(-\frac{J''}{J}) \nonumber \\
&\leq 8\pi (1+ \epsilon)+ 4\pi \mathrm{V}_{M^3}(- 2- \epsilon). \nonumber 
\end{align}

Let $\epsilon\rightarrow 0$ in the above, we obtain
\begin{align}
\varlimsup_{s\rightarrow\infty}\frac{1}{s}\int_{B_p(s)}R\leq 8\pi (1- \mathrm{V}_{M^3}). \label{upper limit has upper bound}
\end{align}

From (\ref{lower limit of global integral}) and (\ref{upper limit has upper bound}), we get
\begin{align}
\lim_{r\rightarrow\infty} \frac{\int_{B_p(r)}R}{r}= 8\pi (1- \mathrm{V}_{M^3}). \nonumber 
\end{align} 

Because $\displaystyle\lim_{r\rightarrow\infty} \frac{\int_{B_p(r)}R}{r}$ exists and is finite, it is easy to see that
\begin{align}
\lim_{r\rightarrow\infty} \frac{\int_{B_q(r)}R}{r}= \lim_{r\rightarrow\infty} \frac{\int_{B_p(r)}R}{r}= 8\pi (1- \mathrm{V}_{M^3}), \quad \quad \quad \forall q\in M^3. \nonumber 
\end{align}
}
\qed

\begin{cor}\label{cor integ of Rc alogn the ray and volume of slices}
{If $(M^3, g)$ is a Riemannian manifold with a pole $p$ and $Rc\geq 0$, then
\begin{align}
\lim_{s\rightarrow\infty} \frac{\mathscr{A}'(s)}{s}= 8\pi\cdot \mathrm{V}_{M^3}. \nonumber 
\end{align}
}
\end{cor}

\pf
{From (\ref{integ equa pf scalar curv}), using Lemma \ref{lem integ of Rc alogn the ray} and Theorem \ref{thm main result}, we have 
\begin{align}
\lim_{s\rightarrow\infty} \frac{\mathscr{A}'(s)}{s}&= \lim_{s\rightarrow\infty} \int_{\rho^{-1}[a, s]}\frac{J''}{s\cdot J}= 8\pi+ \lim_{s\rightarrow\infty}s^{-1}\Big\{\int_{\rho^{-1}[a, s]} Rc(\nabla \rho)- R(M^3) \Big\}  \nonumber \\
&= 8\pi+ \lim_{s\rightarrow\infty}s^{-1}\Big\{\int_{B_p(s)} Rc(\nabla \rho)- R(M^3) \Big\}= 8\pi\cdot \mathrm{V}_{M^3}. \nonumber
\end{align}
}
\qed

\begin{theorem}
{If $(M^3, g)$ is a Riemannian manifold with a pole and Ricci pinched, then $(M^3, g)$ is flat.
}
\end{theorem}

\pf
{From Ricci pinched assumption, there is $\epsilon> 0$ such that $Rc\geq \epsilon\cdot R\cdot g$. Now by Lemma \ref{lem integ of Rc alogn the ray} and Theorem \ref{thm main result}, we get
\begin{align}
0= \lim_{r\rightarrow\infty} \frac{\int_{B_p(r)}Rc(\nabla \rho)}{r}\geq \epsilon\cdot \lim_{r\rightarrow\infty} \frac{\int_{B_p(r)}R}{r}= \epsilon\cdot 8\pi (1- \mathrm{V}_{M^3}) . \nonumber 
\end{align}

Therefore, we obtain that $\mathrm{V}_{M^3}= 1$, which implies that $M^3$ is isometric to $\mathbb{R}^3$ by the rigidity part of Bishop-Gromov's Volume Comparison Theorem.
}
\qed

\end{document}